\documentclass{article}
\usepackage{amssymb,amsmath,amsthm,amscd,verbatim}

\usepackage{hyperref}

\usepackage[shortalphabetic]{amsrefs}

\newtheorem{thm}{Theorem}[section]
\newtheorem{pro}[thm]{Proposition}
\newtheorem{cor}[thm]{Corollary}
\newtheorem{lem}[thm]{Lemma}
\newtheorem{df}[thm]{Definition}
\newtheorem{rem}[thm]{Remark}

\renewcommand{\iff}{\leftrightarrow}

\newcommand{\mc}{\mathcal}

\newcommand{\mf}{\mathfrak}

\newcommand{\restrict}{\upharpoonright}

\renewcommand{\and}{\,\&\,}

\newcommand{\bo}{\mathbf}

\begin{document}

\title{A rigid cone in the truth-table degrees with jump}

\author{Bj\o rn Kjos-Hanssen}

\maketitle

\begin{abstract}
The automorphism group of the truth-table degrees with order and jump is
fixed on the set of degrees above the fourth jump, $ \mathbf 0^{(4)}$.
\end{abstract}

\tableofcontents
\newpage
\section{Introduction}

 A \emph{cone} in a partial order $(D,\le)$ is a set of the form $D(\ge a):=\{x\in D:x\ge a\}$ for some $a\in D$. A subset of $S$ of $D$ is \emph{rigid} if it is fixed under the action of the automorphism group $\text{Aut}(D,\le)$, i.e., for each $x\in S$ and each $\pi\in\text{Aut}(D,\le)$, $\pi(x)=x$. We will also be interested in the case of structures $(D,\le,\mf j)$ where $\mf j$ is a unary function on $D$. In that case, rigidity of $S\subseteq D$ is defined with respect to $\text{Aut}(D,\le,\mf j)$ rather than $\text{Aut}(D,\le)$. 

 It is not known whether the structure of the Turing degrees is rigid, but it is known \cite{JS77} that the structure of the Turing degrees with jump contains a rigid cone. This is shown by applying a jump inversion theorem and results on initial segments. Here we show that also the structure of truth-table degrees with jump $(\mc D_{tt},\le,\mf j)$ contains a rigid cone. For definitions relating to initial segments we refer the reader to the author's doctoral dissertation \cite{K02}, survey article \cite{K03}, and forthcoming article \cite{K11}.

Our main result is that each automorphism of the truth-table degrees with jump is equal to the identity on the cone above $\mathbf 0^{(4)}$. This contrasts with the results of Anderson \cite{Anderson} that each automorphism of the truth-table degrees (not necessarily jump invariant) is equal to the identity on some cone, and each automorphism that preserves $\mathbf 0^{(3)}$ and $\mathbf 0^{(5)}$ is equal to the identity on the cone above $\mathbf 0^{(5)}$. It is still open whether non-trivial automorphisms of these structures exist at all.

\section{Steps of the proof}

In this section we describe the global structure of the proof of our main theorem \ref{3.8}; further recursion-theoretic and lattice-theoretic details are given in the subsequent sections.

\begin{df}\label{3.1} 
In the $tt$-degrees we denote the order by $\leq$. If $\mathbf{x}$, $\mathbf{y}$ are $tt$-degrees, we say that $\mathbf{x}\equiv_T \mathbf{y}$ if for some $X\in \mathbf{x}$ and $Y\in \mathbf{y}$, we have $X\equiv_T Y$.
\end{df}

The following theorem is due to Mohrherr \cite{Mohrherr}.

\begin{thm}\label{3.2} 
Let $n\geq 1$ and $\mathbf{a}\geq \mathbf 0^{(n)}$. Then for some $\mathbf{b}$, $\mathbf{a}=\mathbf{b}^{(n)}$.
\end{thm}

\begin{pro}\label{3.3} 
For each $\mathbf{g}$, $[\bo 0, \mathbf{g}]$ is $\Sigma_{3}^{0}(\bo g)$-presentable.
\end{pro}
\begin{proof}
An analysis of the definition of $tt$-reducibility.
\end{proof}

\begin{cor}\label{3.4} 
Each upper semilattices with least and greatest element that can be realized as initial segments $[\bo 0,\mathbf{g}]$ with $\mathbf{g}^{(2)}\leq \mathbf{y}^{(3)}$ is $\Sigma_{4}^{0}(\bo y)$-presentable.
\end{cor}

\begin{thm}\label{3.5} 
For any $y$, the upper semilattices with least and greatest element that can be realized as initial segments $[\bo 0,\mathbf{g}]$ with $\mathbf{g}^{(2)}\leq \mathbf{y}^{(3)}$ are exactly the $\Sigma_{4}^{0}(\mathbf{y})$-presentable ones.
\end{thm}
\begin{proof}
By Corollary \ref{3.4} and Theorem \ref{1.4}.
\end{proof}

\begin{thm}\label{3.6} 
Let $\pi$ be an automorphism of the truth-table degrees with jump and let $\mathbf{x}\geq \mathbf 0^{(3)}$. Then $\pi(\mathbf{x})\equiv_T {\mathbf{x}}$.
\end{thm}
\begin{proof}
By Theorem \ref{3.2} there is a $\mathbf{y}$ such that $\mathbf{x}=\mathbf{y}^{(3)}$. The initial segments $[\bo 0,\bo y']$ and $[\bo 0,\pi(\bo y')]$ are jump-isomorphic via $\pi$, so by Theorem \ref{3.5}, the $\Sigma^0_4(\bo y)$- and $\Sigma^0_4(\pi(\bo y))$-presentable bounded usls coincide. Hence by Proposition \ref{2.7}, 
$$
\pi(\mathbf{y})^{(3)}\equiv_T \mathbf{y}^{(3)}
$$
and so
$$
\pi(\mathbf{x})=\pi(\mathbf{y}^{(3)})=\pi(\mathbf{y})^{(3)}\equiv_T \mathbf{y}^{(3)}=\mathbf{x}.
$$
\end{proof}

\begin{lem}\label{3.7} 
$\mathbf{a}\equiv_T{\mathbf{b}}\Rightarrow \mathbf{a}'=\mathbf{b}^{\prime}$.
\end{lem}

\begin{thm}\label{3.8} 
Let $\pi$ be an automorphism of the truth-table degrees with jump and let $\mathbf{x}\geq  \mathbf 0^{(4)}$. Then $\pi(\mathbf{x})=\mathbf{x}$.
\end{thm}
\begin{proof}
By Theorem \ref{3.2}, there is a $\mathbf{y}$ such that $\mathbf{x}=\mathbf{y}^{(4)}$. Let $\mathbf{z}=\mathbf{y}^{(3)}$, so $\mathbf{x}=\mathbf{z}^{\prime}$ and $\mathbf{z}\geq  \mathbf 0^{(3)}$. By Theorem \ref{3.6}, $\pi(\mathbf{z})\equiv_T {\mathbf{z}}$ and by Lemma \ref{3.7}, $\mathbf{a}\equiv_T{\mathbf{b}}\Rightarrow \mathbf{a}^{\prime}=\mathbf{b}^{\prime}$. Hence
$$
\pi(\mathbf{x})=\pi(\mathbf{z}^{\prime})=\pi(\mathbf{z})^{\prime}=\mathbf{z}^{\prime}=\mathbf{x}.
$$
\end{proof}

\section{Mal'tsev homogeneous lattice tables}

If $(L,\leq)$ is a partial order (transitive, reflexive, antisymmetric relation) such that greatest lower bounds $\alpha\wedge\beta$ of all $\alpha,\beta\in L$ exist then $(L, \le, \wedge)$ is called a lower semilattice; if least upper bounds $\alpha\vee\beta$ of all pairs $\alpha,\beta\in L$ exist, then $(L, \leq, \vee)$ is called an upper semilattice (usl). If $L$ is both an lower semilattice and an upper semilattice then $L$ is a lattice. $L$ is called bounded if there exist elements $0,1\in L$ such that for all $\alpha\in L$, $0\leq\alpha\leq 1$. In particular every finite lattice is bounded. If $L$ has more than one element (so in the bounded case, $0\neq 1$) then we say that $L$ is nontrivial. A unary algebra is a collection of functions $f$ : $X\rightarrow X$ on a set $X$, closed under composition. The partition lattice Part($X$) on a set $X$ consists of all equivalence relations (considered as sets of ordered pairs) on $X$, ordered by inclusion. We will be interested in the case where $X$ is finite or countably infinite.

 A lattice table (see \cite{Lerman}) $\Theta$ is (1) a set $X$ together with (2) a finite set of equivalence relations $\alpha_{1},\ldots,\alpha_{n}$ on $X$, and (3) an order $\leq$ given by $\alpha_{i}\leq\alpha_{j}\leftrightarrow \alpha_{i}\supseteq\alpha_{j}$ (reverse inclusion of sets of ordered pairs), such that $\{\alpha_{1},\ldots,\alpha_{n}\}$ ordered by \emph{inclusion} is a 0-1 sublattice of Part(X). We write $\widehat{\Theta}=\{\alpha_{1},\ldots,\alpha_{n}\}$. We think of $\Theta$ as equal to $X$, but endowed with additional structure. So $ x\in\Theta$ means $x\in X$, etc. but for emphasis we may write $|\Theta|$ for $X$. Note that $\Theta$ is determined by $\widehat{\Theta}$.

 Elements of $|\Theta|$ are denoted by lower-case Roman letters such as $u$, $v$, $w$, $x$, $y$, $z$, and elements of semilattices in general and $\widehat{\Theta}$ in particular by lower-case Greek letters such as $\alpha$, $\beta$, $\gamma$.

 If $\alpha\in\widehat{\Theta}$ and $(x, y)\in\alpha$ then we write $x\sim_{\alpha}y$. If $\Theta$ is a lattice table then an endomorphism of $\Theta$ is a map from $\Theta$ to $\Theta$ preserving all equivalence relations in $\widehat{\Theta}$. That is, $(\forall x, y\in\Theta)(\forall\alpha\in\widehat{\Theta})(x\sim_{\alpha}y\rightarrow f(x)\sim_{\alpha}f(y))$. End $\Theta$ denotes the unary algebra consisting of all endomorphisms of $\Theta$.

 $C_{\Theta}(x, y)$ denotes the principal equivalence relation in $\Theta$ generated by $(x, y)$, i.e.
$$
C_{\Theta}(x, y)=\cap\{\alpha\in\widehat{\Theta} : (x, y)\in\alpha\}.
$$
We define $\text{End}_{\Theta}(x, y)$ to be the the principal congruence relation in $\Theta$ generated by $(x, y)$, i.e. the equivalence relation generated by all pairs $(f(x), f(y))$ for $f\in \text{End }\Theta$.

\begin{lem}\label{4.1} 
$\text{End}_{\Theta}(x, y)\subseteq C_{\Theta}(x, y)$.
\end{lem}
\begin{proof}
If $(u, v)\in \text{End}_{\Theta}(x, y)$ then $(c, d)$ is in the transitive closure of $$\{(f(x), f(y))\mid f\in\text{ End }\Theta\},$$ so it suffices to show each such $(f(x), f(y))\in C_{\Theta}(x, y)$. For this it suffices to show $(f(x), f(y))\in\alpha$ provided that $(x, y)\in\alpha$ for $\alpha\in\widehat\Theta$; this holds since $f\in \text{End } \Theta$.
\end{proof}

\begin{df}\label{4.2} 
Let $\Theta$ be a lattice table. We say that $\Theta$ is \emph{Mal'tsev homogeneous} if for all $x,  y\in\Theta, C_{\Theta}(x, y)\subseteq \text{End}_{\Theta}(x, y)$ (so by Lemma \ref{4.1}, $C_{\Theta}(x, y)=\text{End}_{\Theta}(x, y))$.
\end{df}

The following Proposition can readily be proved:

\begin{pro}\label{4.3} 
$\Theta$ is Mal'tsev homogeneous iff for all $x, y, u,  v\in\Theta$ satisfying
$$
(\forall\alpha\in\widehat{\Theta})(x\sim_{\alpha}y\rightarrow u\sim_{\alpha}v),
$$
there exist $n\in\omega=\{0,1,2,\ldots\}, z_{1},\ldots,  z_{n}\in\Theta$ and $f_{0},\ldots,  f_{n}\in$ End $\Theta$ such that
$$
(\forall i\leq n)(\{f_{i}(x), f_{i}(y)\}=\{z_{i}, z_{i+1}\})
$$
where $z_{0}=u$ and $z_{n+1}=v$. 
\end{pro}

The $z_{i}$ are called \emph{homogeneity interpolants}.

 This notion of homogeneity is more general (weaker) than those considered in \cite{Lerman}.

 Note that if $\alpha\wedge\beta=\gamma$ in $\widehat{\Theta}$ then $\alpha$ and $\beta$ generate $\gamma$. That is, if $x\sim_{\gamma}y$ then there exist \emph{meet interpolants} $z_{1},\ldots, z_{n}$ for $x,y$ such that $x\sim_{\alpha}z_{1}\sim_{\beta} z_{2}\cdots\sim_{\alpha} z_{n}\sim_{\beta} y$.

\begin{df}\label{4.4} 
If $\Theta$ is a lattice table and $Y\subseteq|\Theta|$, then for each $\alpha\in\widehat{\Theta}$, $\alpha\restrict Y=\{(x, y)\in Y\times Y : (x, y)\in\alpha\}$. Let $\widehat{\Theta}\restrict Y=\{\alpha\restrict Y :\alpha\in\widehat{\Theta}\}$.

 If $\Theta_{0}$ and $\Theta_{1}$ are lattice tables, then we say that $\Theta_{0}\subseteq\Theta_{1}$ if $|\Theta_{0}|\subseteq|\Theta_{1}|$ and $\widehat{\Theta}_{1}\restrict |\Theta_{0}|=\widehat{\Theta}_{0}$. Note that if $\Theta_{0}\subseteq\Theta_{1}$ then $\widehat{\Theta}_{0}$ and $\widehat{\Theta}_{1}$ are isomorphic (nontrivial, finite) lattices.

 If $\Theta_{n},  n\in\omega$ are lattice tables such that $\Theta_{n}\subseteq\Theta_{n+1}$ for each $n$, then $ \bigcup_{n\in\omega}\Theta_{n}$ is the lattice table $\Theta$ such that $| \Theta|=\bigcup_{n\in\omega}|\Theta_{n}|$ and $\widehat{\Theta}\restrict |\Theta_{n}|=\widehat{\Theta}_{n}$ for each $n$. In particular $\Theta_{n}\subseteq\Theta$ and $\widehat{\Theta}_{n}$ and $\widehat{\Theta}$ are isomorphic lattices for each $n$.
\end{df}

\begin{df}\label{4.5} 
$\Theta$ is a sequential lattice table if there exist $\Theta_{n},  n\in\omega$, such that $ \Theta=\bigcup_{n\in\omega}\Theta_{n}$, and

\begin{enumerate}
\item[(1)] each $\Theta_{n}$ is a $(0,1,\vee)$-substructure of Part $(|\Theta_{n}|)$ ( $\Theta_{n}$ is an \emph{usl table}),

\item[(2)] $\Theta$ is a lattice table, and

\item[(3)] for each $n$, meet interpolants for elements of $\Theta_{n}$ exist in $\Theta_{n+1}$.
\end{enumerate}

$\Theta$ is a sequential Mal'tsev homogeneous lattice table if in addition

\begin{enumerate}
\item[(4)] $\Theta$ is Mal'tsev homogeneous, with homogeneity interpolants for elements of $\Theta_{n}$ appearing in $\Theta_{n+1}$ (compare VII.1.1, 1.3 of \cite{Lerman}).
\end{enumerate}
\end{df}

\begin{df}[Direct limit]\label{4.6} 
Let a sequence $(L^{i},\varphi_{i})_{i\in\omega}$ be given, where each $L^{i}$ is a finite lattice, $\varphi_{i}$ : $L^{i}\rightarrow L^{i+1}$ is a $(0,1,\vee)$ homomorphism, and $L^{i}\cap L^{j}=\emptyset$ for $i\neq j$.

 Let $L^{\prime}= \bigcup_{i\in\omega}L_{i}$ as a set. Let $\approx$ be the equivalence relation on $L^{\prime}$ generated by $a\approx\varphi_{i}(a)$ for $a\in L^{i}$. Then $ L=L^{\prime}/\approx$ is an upper semilattice called the direct limit of the sequence $(L^{i},\varphi_{i})_{i\in\omega}$.
\end{df}

\begin{df}\label{4.7} 
Fix finite lattices $L^{0}, L^{1}$ and a $(0,1, \vee)$ homomorphism $\varphi$ : $L^{0}\rightarrow\varphi(L^{0})\subseteq L^{1}$, and lattice tables $\Theta^{0},\Theta^{1}$. Suppose $\Psi^{i}$ : $L^{i}\rightarrow\widehat{\Theta}^{i}, i=0,1$, are isomorphisms. For $\alpha\in L^{i}$, we write $\sim_{\alpha}$ for $\sim_{\Psi^{i}\alpha}$.

 We say that \emph{$\Theta^{1}$ embeds in $\Theta^{0}$ with respect to $\varphi$ and $\Psi_{0},\Psi_{1}$} if there is a function $\Theta(\varphi)$ : $\Theta^{1}\rightarrow\Theta^{0}$ such that for all $x, y\in\Theta^{1}$, and all $\alpha\in L^{0}$,
$$
x\sim_{\varphi\alpha}y\Leftrightarrow\Theta(\varphi)(x)\sim_{\alpha}\Theta(\varphi)(y).
$$
\end{df}

 In the characterization of intervals $[\bo 0, \mathbf{g}]$ for $ \mathbf{g}< \mathbf 0^{\prime}$, Definition \ref{4.7} plays a key role which we will now describe.

 Suppose a bounded countable upper semilattice $L$ is given such that the ordering $\leq$ of $L$ is computably enumerable but not necessarily computable. That is, there is a computable sequence that consists of all pairs $(\alpha,\beta)$ such that $\alpha\leq\beta$, but if a given pair $(\alpha,\beta)$ does not appear anywhere in the list then this cannot be determined effectively.

 For reasons whose explanation would take us too far afield (but see \cite{K03}), we need a computable sequence of lattice tables $\Theta^{0},\Theta^{1},\ldots$ such that $\widehat\Theta^{s}$ is isomorphic to our approximation to $L$ at stage $s$. (We will start with a sequence $(L^{i},\varphi_{i})_{i\in\omega}$ having $L$ as direct limit, and our approximation to $L$ at stage $s$ will be $L^{s}.$) Suppose we discover at stage $s+1$ that $\alpha\leq\beta$, whereas at stage $s$ we knew that $\beta\leq\alpha$ but thought that $\alpha\not\leq\beta$. Further suppose that we cannot ignore what was done using $\Theta^{s}$ at stage $s$, but we can let $\Theta^{s+1}$ be a subset of $\Theta^{s}$. If $\Theta^{s+1}$ embeds into $\Theta^{s}$ with respect to $\varphi$ (the homomorphism mapping our approximation to $L$ at stage $s$ to our approximation to $L$ at stage $s+1$) then by thinning $\Theta^{0}$ to $\Theta(\varphi)\Theta^{1}$, we eliminate all elements $x, y$ that are witnesses to the fact that $\alpha\neq\beta$. This allows us to identify $\alpha$ and $\beta$, even though so far we have been working under the assumption that $\alpha\neq\beta$.

We mention for the reader who is a computability theorist that in the characterization of lattices isomorphic to $[\bo 0, \mathbf{g}]$ for $ \mathbf{g}< \mathbf 0^{\prime}$, the ordering of $L$ is computably enumerable only relative to the Turing degree $\bo 0''$, and the ``stages'' above are really levels of a priority tree, the true path of which it requires $\bo 0''$ to identify.

 The full result needed for the application to the Turing degrees is contained in Theorem \ref{4.8} and Proposition \ref{4.9}.

\begin{thm}\label{4.8} 
Let $L$ be a bounded countable nontrivial usl and let $(L^{i},\varphi_{i})_{i\in\omega}$ be any system of nontrivial finite lattices having $L$ as direct limit in the sense of Definition \ref{4.6}. Then there exists

\begin{enumerate}
\item[1.] a function $ h:\omega\rightarrow\omega$,

\item[2.] a double sequence of finite lattice tables $(\Theta_{j}^{i})_{i\in\omega,j\geq h(i)}$ with $\Theta_{j}^{i}\subseteq\Theta_{j+1}^{i}$
 for each $ i\in\omega, j\geq h(i)$, and

\item[3.] for each $ i\in\omega$ an increasing function $m_{i}$ : $\omega\rightarrow\omega$ with $m_{i}(0)=h(i)$, such that
\begin{enumerate}
\item[1.] letting $ \Theta^{i}=\bigcup_{j\in\omega}\Theta_{j}^{i}$, we have $|\Theta^{i}|\supseteq|\Theta^{i+1}|$ for each $ i\in\omega$,

\item[2.] for each $i, j\geq h(i)$ and $k$ such that $m_{i}(j)\leq k<m_{i}(j+1)$, we have
$$
\Theta_{k}^{i}=\Theta_{m_{i}(j)}^{i},
$$
\item[3.] for each $ i\in\omega, (\Theta_{m_i(j)}^{i})_{j\in\omega}$ is a sequential Mal'tsev homogeneous lattice table,

\item[4.] for each $i\in\omega,\widehat{\Theta}_{i}$ is isomorphic to $L^{i}$, and

\item[5.] there exist isomorphisms $\Psi_{i}$ : $L^{i}\rightarrow\widehat{\Theta}_{i}$ such that $\Theta^{i+1}$ embeds in $\Theta^{i}$ with respect to $\varphi_{i}$ and $\Psi_{i},\Psi_{i+1}$, and the embedding is the identity map. In
 other words, for all $x, y\in\Theta^{i+1}$ and $\alpha\in L^{i}$, we have
$$
x\sim_{\Psi_{i}\alpha} y\leftrightarrow x\sim_{\Psi_{i+1}\varphi_{i}\alpha} y.
$$

\end{enumerate}
\end{enumerate}

\end{thm}

 The essential property in Theorem \ref{4.8}, and the one that goes beyond those of \cite{LL76}, is (5). The following Proposition can be proved by inspecting the proof of Theorem \ref{4.8}.

\begin{pro}[Computability-theoretic properties] \label{4.9}
Let $\bo a$ be a Turing degree and let $L$ be a $\Sigma_{1}^{0}(\bo a)$-presentable usl, as in \cite{K03}. Then in Theorem \ref{4.8}, we may assume that $h$ is $\bo a$-computable; the array $\{\Theta_{j}^{i}\mid j\geq h(i)\}$ is $\bo a$-computable; each $m_{i}$ is computable; for each $ i<\omega, (\Theta_{m_{i}(j)}^{i})_{j\in\omega}$ is a computable sequence; each $\Theta^{i}$ is computable; and there is a computable function taking $L^{0},\ldots, L^{i}$ to $\Theta^{i}$.
\end{pro}

 We now begin the development that will lead to a proof of Theorem \ref{4.8}.

 If $A$ is a unary algebra then Con $A$ denotes the congruence lattice of $A$, i.e. the lattice of all equivalence relations $E$ on $X$ preserved by all $f\in A$, ordered by inclusion.

 The following observation can be traced back to Mal'tsev \cites{UA, Maltsev54, Maltsev63}.

\begin{pro}\label{4.10} 
For any unary algebra $A$, the dual of Con $A$ is a Mal'tsev homogeneous lattice table.
\end{pro}
\begin{proof}
Suppose $A$ is a unary algebra on a set $X$. Let $\Theta$ be the lattice table such that $\widehat\Theta$ is the dual of Con $A$. Since Con $A$ is a 0-1 sublattice of Part($A$), $\Theta$ is a lattice table.

 If $f$ is an operation in $A$ and $\alpha\in\widehat{\Theta}$ then $\alpha$ is a congruence relation on $A$ and hence $\forall x, y(x\sim_{\alpha}y\rightarrow f(x)\sim_{\alpha}f(y))$, which means that $f\in \text{End } \Theta$. So
$$
A\subseteq \text{End }\Theta.
$$
Clearly for any unary algebras $A, B$ on the same underlying set, we have $ A\subseteq  B\Rightarrow \text{Con }A\supseteq\text{Con }B$. Hence Con End $\Theta\subseteq$ Con $A=\widehat{\Theta}$.

 If $u$, $v$, $x$, $y\in X$, $ f\in$ End $\Theta$ and $(x, y)\in \text{End}_{\Theta}( u, v)$ then there exist $z_{1},\ldots, z_{k}$ such that $(z_{i}, z_{i+1})=(g_{i}(u), g_{i}(v))$ for $g_{i}\in \text{End } \Theta$ with $z_{1}=x, z_{k}= y$, hence letting $w_{i}=f(z_{i})$ and $h_{i}=f\circ g_{i}$ we have $(w_{i}, w_{i+1})=(h_{i}(u), h_{i}(v))\in \text{End}_{\Theta}(u, v)$, $w_{1}=f(x), w_{k}=f(y)$, and so $(f(x), f(y))\in \text{End}_{\Theta}(u, v)$. Hence we have shown $\text{End}_{\Theta}(u, v)\in$ Con End $\Theta\subseteq\Theta$. Since End $\Theta$ contains the identity map, $\text{End}_{\Theta} (u, v)$ contains $(u, v)$. Hence $\text{End}_{\Theta}(u, v)$ is in $\Theta$ and contains $(u, v)$, so it contains $C_{\Theta}(u, v)$. So $\Theta$ is Mal'tsev homogeneous. 
\end{proof}

We recall the construction of \cite{Pudlak76}.

\begin{df}\label{4.11} 
Let $L$ be a nontrivial lattice. $\mathcal{A}=(A, r, h)$ is called an $L-\{1\}$-\emph{colored} graph if $A$ is a set, $r$ is a set of size-two subsets of $A$, i.e. $(A, r)$ is an undirected graph without loops, and $h:r\rightarrow L-\{1\}$ is a mapping of the set $r$ of the edges of the graph into $L-\{1\}$.

 The map $e$ : $L\rightarrow \text{Part}(A)$ is defined by: for $\alpha\in L, e(\alpha)$ is the equivalence relation on $A$ generated by identifying points $x, y$ if there is a path from $x$ to $y$ in the graph consisting of edges all of which have color $\geq\alpha$. In this case we say that $x, y$ are connected with color $\geq\alpha$.
\end{df}

\begin{df}[$\alpha$-cells]\label{4.12}
Let $L$ be a nontrivial lattice and let $\alpha\in L-\{1\}$. An $\alpha$-cell $\mathcal{B}_{\alpha}=(B_{\alpha}, s_{\alpha}, k_{\alpha})$ is an $L-\{1\}$-colored graph consisting of (1) a base edge $\{x, y\}$ colored $\alpha$, and (2) for each pair $(\alpha_{1},\alpha_{2})$ of elements of $L$ such that $\alpha_{1}\wedge\alpha_{2}\leq\alpha$, a chain of edges $\{x, u_{1}\},\{u_{1}, u_{2}\},\{u_{2}, u_{3}\},\{u_{3}, y\}$, colored $\alpha_{1},\alpha_{2},\alpha_{1},\alpha_{2}$, respectively. Here $x, y, u_{1}, u_{2}, u_{3}$ are distinct elements of $B_{\alpha}$. The base edge and chain of edges corresponding to a particular inequality $\alpha_{1}\wedge\alpha_{2}\leq \alpha$ is referred to as a pentagon. So an $\alpha$-cell consists of several pentagons, intersecting only in a common base edge.
\end{df}

\begin{df}[Pudl\'{a}k graphs]\label{4.13} 
Let $L$ be a nontrivial lattice. The Pudl\'{a}k graph \cite{Pudlak76} of $L$ is an $L-\{1\}$ colored graph $\mathcal{A}^{P}$, defined as follows.

\begin{enumerate}
\item[1.] $\mathcal{A}_{0}^{P}$ consists of a single edge colored by $0\in L$. (In fact, how we choose to
 color this one edge has no impact on later proofs.)

\item[2.] $\mathcal{A}_{n+1}^{P}$ contains $\mathcal{A}_{n}^{P}$ as a subgraph and is obtained by attaching to each edge of $\mathcal{A}_{n}^{P}$ of any color $\alpha$ an $\alpha$-cell.

\item[3.] $ \mathcal{A}^{P}=\bigcup_{n\in\omega}\mathcal{A}_{n}^{P}$.
\end{enumerate}

We will use the following modification, which contains infinitely many copies of each edge in Pudl\'{a}k's graph.

\begin{enumerate}
\item[1.] $\mathcal{A}_{0}^{(i)}=\mathcal{A}_{0}^{P}$, for each $ i\in\omega$.

\item[2.] $\mathcal{A}_{j}^{(i)}$ is obtained by attaching to each edge of $\mathcal{A}_{j-1}^{(i)}$ of any color $\alpha, i$ many
$\alpha$-cells.

\item[3.] $\mathcal{A}_{j}=\mathcal{A}_{j}^{(j)}$.

\item[4.] $ \mathcal{A}=\bigcup_{n\in\omega}\mathcal A_{n}=\mathcal{A}(L)$ is called the homogenized Pudl\'{a}k graph of $L$.
\end{enumerate}

The underlying set of $\mathcal A_{n}$ is denoted by $A_{n}$.

Let $\Theta=\Theta(L)$ be the lattice table with $|\Theta|=A$, and $\widehat{\Theta}=\{e(\alpha) :\alpha\in L\}$.

Note that by definition of $\Theta$ being a lattice table, $\widehat{\Theta}$ is ordered by reverse inclusion. Similarly let $\Theta_{n}$ be the lattice table with $|\Theta_{n}|=A_{n}$, $$\widehat{\Theta}_{n}=\{e(\alpha)\restrict \Theta_{n}\mid\alpha\in L\}.$$
\end{df}

\begin{lem}\label{4.14} 
Let $B_{0}\subseteq B_{1}\subseteq A$, where $A$ is the underlying set of $\Theta(L)$.

 For $i=0,1$, let $\Xi_i$ be the usl table whose underlying set is $B_{i}$, and whose equivalence relations are computed using graph points belonging to $B_{i}$ only.

 Then $\Xi_0\subseteq\Xi_1$ in the sense of Definition \ref{4.4}.
\end{lem}
\begin{proof}
We have to show that if $x\sim_{\alpha}y$ holds in $\Xi_1$ then $x\sim_{\alpha}y$ holds in $\Xi_0$.
 The only way this could fail is if there is a path of edges between $x$ and $y$ leading out of $\Xi_0$ and then back in. We may assume that the path does not leave and re-enter $\Xi_0$ via the same node. So it suffices to show that any path that goes around a pentagon not contained in $\Xi_0$ but whose base is in $\Xi_0$ can be shortened to one contained in $\Xi_0$ with no loss of equivalence. Since the pentagons represent inequalities $\alpha\wedge\beta\leq\gamma$, any path $x, u_{1}, u_{2}, u_{3}, y$ going around the pentagon in $\Xi_1- \Xi_0$ may be replaced by the edge $x, y$ cutting across which has equal or greater color, i.e. with no loss of equivalence. 
\end{proof}

\begin{lem}\label{4.15} 
$\Theta_{n}\subseteq\Theta_{n+1}$ for each $n$, so $ \Theta=\bigcup_{n\in\omega}\Theta_{n}$.
\end{lem}
\begin{proof} 
Let $\Xi_i=\Theta_{n+i}$ for $i=0,1$ and apply Lemma \ref{4.14}. 
\end{proof}

\begin{thm}\label{4.16} 
Let $L$ be a nontrivial finite lattice. $L$ is dual isomorphic to the congruence lattice of End $\Theta(L)$. In fact, $e:L\rightarrow\widehat{\Theta}(L)$ is an isomorphism, and $\Theta(L)=\text{Con}$ End $\Theta(L)$.
\end{thm}
\begin{proof}
Pudl\'{a}k \cite{Pudlak76} assumes that $L$ is an algebraic lattice \cite{GLT}, defines a certain algebra $S \subseteq$ End $\Theta^{P}(L)$, and shows that $e$ : $L\rightarrow\widehat{\Theta}^{P}(L)$ is an isomorphism, and $\widehat{\Theta}^{P}(L)=$ Con $S$. Now trivially $\widehat{\Theta}^{P}(L)\subseteq$ Con End $\Theta^{P}(L)$ holds, and $S \subseteq$ End $\Theta^{P}(L)$ implies Con End $\Theta^{P}(L)\subseteq$ Con $S$. So we have $\widehat{\Theta}^{P}(L)=\text{Con}$ End $\Theta^{P}(L)$.

 In fact Pudl\'{a}k's proof works for our graph $\Theta$ as well, i.e. it shows that $e:L\rightarrow\widehat{\Theta}(L)$ is an isomorphism, and $\widehat{\Theta}(L)=\text{Con}$ End $\Theta(L)$.

 Now let $L$ be a finite lattice. Since every finite lattice is algebraic, $L$ is an algebraic lattice. Hence $e$ : $L\rightarrow\widehat{\Theta}(L)$ is an isomorphism and $\widehat{\Theta}(L) =$ Con End $\Theta(L)$.
\end{proof}

\begin{lem}\label{4.17} 
The sequence $\Theta_{n}(L),  n\in\omega$ has a subsequence which is a com- putable Mal'tsev homogeneous sequential lattice table.
\end{lem}
\begin{proof}
Since $\Theta$ is a congruence lattice, $\Theta$ is a Mal'tsev homogeneous lattice table. Hence as $ \Theta=\bigcup_{n\in\omega}\Theta_{n}$, a subsequence of $\Theta_{n},  n\in\omega$ will be a sequential Mal'tsev homogeneous lattice table. The sequence is computable since to compute an equivalence relation on elements of $\Theta_{n}$, it is sufficient to consider paths in $\Theta_{n}$, since $\Theta_{n}\subseteq\Theta_{n+1}$ by Lemma \ref{4.15}.
\end{proof}

 From now on we will assume that in fact $\Theta_{n},  n\in\omega$ is itself the subsequence from Lemma \ref{4.17}.
 Fix nontrivial finite lattices $L^{0}, L^{1}$ and a $(0,1, \vee)$-isomorphisms $\varphi$ : $ L^{0}\rightarrow \varphi(L^{0})\subseteq L^{1}$. Let the $\wedge$-isomorphism $\varphi^{*}$ : $L^{1}\rightarrow L^{0}$ be defined by $\varphi^{*}\beta= \bigvee\{\alpha\in L^{0}\mid \varphi(\alpha)\leq\beta\}$. This $\varphi^{*}$ is known as the Galois adjoint of $\varphi$ \cite{Gierz}.

 The map $\varphi^{*}$ has many nice properties; we list the ones we need in the following lemma.

\begin{lem}\label{4.18} 
\begin{enumerate}
\item[1.] \label{4.18.1} $\varphi^{*}$ is a $(\wedge, 1)$-homomorphism.

\item[2.] \label{4.18.2} If $\beta<1$ then $\varphi^{*}\beta<1$.

\item[3.] \label{4.18.3} $\varphi^{*}$ is injective on $\varphi L^{0}$.

\item[4.]  $\alpha\leq\varphi^{*}\beta\leftrightarrow\varphi^{*}\varphi\alpha\leq\varphi^{*}\beta$.\label{4.18.4}
\end{enumerate}
\end{lem}

\begin{proof} 
These all follow easily from the definition of $\varphi^{*}$ and the fact that $\{\alpha\in L^{0}\mid \varphi(\alpha)\leq\beta\}$ is the principal ideal generated by $\varphi^{*}(\beta)$, i.e. $\{\alpha\in L^{0}\mid \alpha \leq\varphi^*(\beta)\}$.
\end{proof}

\begin{lem}\label{4.19} 
Let $\mathfrak{C}(\varphi)\mathcal{A}L^{1}$ be the graph obtained from $\mathcal{A}L^{1}$ by replacing each color $\beta$ by $\varphi^{*}\beta$. Then $\mathfrak{C}(\varphi)\mathcal{A}L^{1}$ is isomorphic to a subgraph of $\mathcal{A}L^{0}$.
\end{lem}
\begin{proof}
Each pentagon of $\mathfrak{C}(\varphi)\mathcal{A}L^{1}$ represents an inequality of the form
$$
\varphi^{*}\beta_{1}\wedge\varphi^{*}\beta_{2}\leq\varphi^{*}\beta,
$$

 for $\beta_{1},\beta_{2},\beta\in L^{1}$ satisfying $\beta_{1}\wedge\beta_{2}\leq\beta$. Then $\varphi^{*}\beta_{1}\wedge\varphi^{*}\beta_{2}=\varphi^{*}(\beta_{1}\wedge\beta_{2})\leq \varphi^{*}\beta$, so the represented inequality $\varphi^{*}\beta_{1}\wedge\varphi^{*}\beta_{2}\leq\varphi^{*}\beta$ holds in $L^{0}$.

 Hence we can obtain an isomorphic copy of $\mathfrak{C}(\varphi)\mathcal{A}L^{1}$ within $\mathcal{A}L^{0}$ by running through the construction of $\mathcal{A}L^{0}$, omitting every pentagon that represents an inequality involving members of $L^{0}-\varphi^{*}L^{1}$, and omitting pentagons for inequalities that are true in $L^{0}$ but not in $L^{1}$. If an edge becomes disconnected from $\mathcal{A}_{0}$ by such omissions then it too is omitted. Since $L^{1}$ may have many more elements than $L^{0}$, we make use of the fact that $\mathcal{A}$ contains infinitely many copies of each edge from Pudl\'{a}k's original graph $\mathcal{A}^{P}$. Since $\varphi^{*}(\beta)=1\rightarrow\beta=1$ by Lemma \ref{4.18}, recoloring of points is never identification of points. 
\end{proof}

\begin{lem}\label{4.20} 
Let $\Theta(\varphi)$ be the isomorphism from Lemma \ref{4.19}, sending $\mathcal{A}L^{1}$ to a subgraph of $\mathcal{A}L^{0}$ isomorphic to $\mathfrak{C}(\varphi)\mathcal{A}L^{1}$.

 Then $\Theta(\varphi)\Theta L^{1}\subseteq\Theta L^{0}$ in the sense of Definition \ref{4.4}.
\end{lem}
\begin{proof}
Let $\Xi_0=\Theta(\varphi)\Theta L^{1}$ and $\Xi_1=\Theta L^{0}$ and apply Lemma \ref{4.14}. 
\end{proof}

\begin{lem}\label{4.21} 
Let $\Psi_{i}$ be the map $e$ of Definition \ref{4.11} for $L=L^{i}, i=0,1$.

 Then $\Theta(L^{1})$ embeds in $\Theta(L^{0})$ with respect to $\varphi$ and $\Psi_{0},\Psi_{1}$.
\end{lem}
\begin{proof}
Let $x, y$ be points in $\Theta L^{1}$, i.e. in $\mathcal{A}L^{1}$, and let $\alpha\in L^{0}$. Then obviously $x\sim_{\varphi\alpha}y\rightarrow\Theta(\varphi)x\sim_{\varphi^{*}\varphi\alpha}\Theta(\varphi)y$. Now suppose $\Theta(\varphi)x\sim_{\varphi^{*}\varphi\alpha}\Theta(\varphi)y$. Then there is a path witnessing this, which by Lemma \ref{4.19} we may assume lies within $\Theta(\varphi)\mathcal{A}L^{1}$. Hence the path has an inverse image path under $\Theta(\varphi)^{-1}$. This is then a path from $x$ to $y$ with colors $\beta$ for all of which $\varphi^{*}\beta\geq\varphi^{*}\varphi\alpha$. But then $\alpha\leq\varphi^{*}\beta$ by Lemma \ref{4.18}(4), and so $\varphi\alpha\leq\beta$, so $x\sim_{\varphi\alpha}y$. So in fact $x\sim_{\varphi\alpha}y\iff \Theta(\varphi)x\sim_{\varphi^*\varphi\alpha}\Theta(\varphi)y$. Colors $\gamma$ of edges in $\Theta(\varphi)\mathcal A L^1$ are all of the form $\varphi^*(\beta)$ for some $\beta$. So $\Theta(\varphi)x\sim_{\varphi^{*}\varphi\alpha}\Theta(\varphi)y$ iff there is a path from $\Theta(\varphi)x$ to $\Theta(\varphi)y$, all edges of which are colored $\gamma\geq\varphi^{*}\varphi\alpha$, or equivalently by Lemma \ref{4.18}(4) (using $\gamma=\varphi^{*}\beta)$, colored $\gamma\geq\alpha$. Hence equivalently $\Theta(\varphi)x\sim_{\alpha}\Theta(\varphi)y$.
\end{proof}

\begin{proof}[Proof of Theorem \ref{4.8}] 
Let $m_{i}(n)$ be the least $m$ such that $\Theta(\varphi_{i})\Theta_{n}^{i+1}\subseteq\Theta_{m}^{i}$. Let $h(i)=m_{i}(0)$, for $ i\in\omega$. Let $\Theta^{0}=\Theta(L^{0})$ and for $i\geq 1$, denoting composition by juxtaposition,
$$
\Theta^{i}=\Theta(\varphi_{0})\cdots\Theta(\varphi_{i-1})\Theta(L^{i}).
$$
Let $\Theta_{k}^{i}=\Theta(\varphi_{0})\cdots\Theta(\varphi_{i-1})\Theta_{j}(L^{i})$ if $k=m_{0}m_{1}\cdots m_{i-1}(j)$ for some $j$; other- wise, let $\Theta_{k}^{i}=\Theta_{k-1}^{i}$. The Theorem now follows easily. 
\end{proof}

\section{Initial segments of the $tt$-degrees}

\begin{lem}\label{1.1}
Suppose for each $e$, $g$ lies on a tree $T_{e}$ which is e-splitting for some $c$ for some tables with the properties of Proposition \ref{4.9}, in the sense of \cite{LL76}. Then $\mathbf{g}$ is hyperimmune-free.
\end{lem}
\begin{proof}
For each $e \in\omega$ there exists $ e^{*}\in\omega$ such that for all stages $s$ and all oracles $g$, if $\{e^{*}\}_{s}^{g}(x)\downarrow$ then $\{e^{*}\}^{g}(x)=\{e\}^{g}(x)$ and $\{e\}_{s}^{g}(y)\downarrow$ for all $y\leq x$. If $g$ lies on $T_{e^{*}}$ then it follows that $\{e\}^{g}$ is total and $\{e^{*}\}^{T(\sigma)}(x)\downarrow$ for each $\sigma$ of length $x+1$. Hence $\{e\}^{g}=\{e^{*}\}^{g}$ is dominated by the recursive function $f(x)= \max\{\{e\}^{T(\sigma)}(x) : |\sigma|=x+1\}$.
\end{proof}

\begin{pro}\label{1.2} 
Let $L$ be a $\Sigma_{4}^{0}(\bo y)$-presentable upper semilattice with least and greatest element.
Then there exist $t$, $i$, $g$ such that
\begin{enumerate}
\item[1.] $ t:\omega \rightarrow 2$ is $0''$-computable,
\item[2.] $i$ is the characteristic function of a set $I$ such that $I\leq_{m}y^{(3)}$,
\item[3.] $g^{(2)}(e)=t(i(0),\ldots, i(e))$ for all $e \in\omega$,
\item[4.] $[\bo 0, \bo g]$ is isomorphic to $L$, and
\item[5.] $\bo g$ is hyperimmune-free
\end{enumerate}
\end{pro}
\begin{proof}
The proof in \cite{LL76} must be modified to employ the lattice tables of Proposition \ref{4.9}.

By Proposition \ref{4.9}, for all $ x$, $y\in\Theta^{k+1}$ and $\alpha\in L^{k}$, we have [identifying the isomorphism between $L^{i}$ and $\widehat\Theta^{i}$ with the identity]
$$
x\sim_{\varphi_{k}\alpha}y\leftrightarrow x\sim_{\alpha}y.
$$
Lemma 4.1 of \cite{LL76} is modified so that $\psi_{T,c}$ is $\psi_{T,\varphi_kc}$. The equivalence
$$
uF_{m(i)}(c)v\leftrightarrow uG_{i}(c)v
$$
now becomes
$$
uF_{m(i)}(c)v\leftrightarrow uG_{i}(\varphi_{k}c)v
$$

 Just as in Lemma \ref{4.1} it is shown that $\psi_{T,c}$ is Turing equivalent to $\psi_{T_{0},c}$, it now follows that $\psi_{T,\varphi_k c}$ is Turing equivalent to $\psi_{T_{0},c}$, which is what we want.

 $i(e)=1$ iff the answers to the $\Pi_{1}^{0}(y^{(2)})$-question about $L^{e}$ is yes.

 By Lemma \ref{1.1}, $\bo g$ is hyperimmune-free.
\end{proof}

\begin{lem}\label{1.3}
If $t$, $i$, $A$, $q$ satisfy
\begin{enumerate}
\item[1.] $t:\omega\rightarrow 2$ is q-computable,
\item[2.] $i$ is the characteristic function of a set $I$ such that $I\leq_{m}q^{\prime}$,
\item[3.] $A(e)=t(i(0),\ldots, i(e))$ for all $ e\in\omega$,
\end{enumerate}
then $A\leq_{tt}q^{\prime}$.
\end{lem}
\begin{proof}
The value of $A(e)$ is determined by the following $e+2$ many yes-or-no questions to $q^{\prime}$:

 Is $ i(0)=0$? $\cdots$ Is $i(e)=0$? and, using the answers to the first $e+1$ many questions: Is $t(i(0),\ldots, i(e))=0$?
\end{proof}

\begin{thm}\label{1.4}
Each $\Sigma_{4}^{0}(\bo y)$-presentable upper semilattice with least and greatest element can be realized as an initial segment $[\bo 0,\bo g]$ with $\bo g^{(2)}\leq \bo y^{(3)}$.
\end{thm}
\begin{proof}
Let $g$ be as in Proposition \ref{1.2}. By Lemma \ref{1.3} with $q=y^{(2)}$ and $A=g^{(2)}$, we have $g^{(2)}\leq_{tt}y^{(3)}$.

 By Proposition \ref{1.2}, $L$ is isomorphic to $[\bo 0, \bo g]_{T}$. Since $\bo g$ is hyperimmune-free, $[\bo 0, \bo g]_{T}=[\bo 0, \bo g]_{tt}$.
\end{proof}

\section{Coding a set into a lattice}

\begin{df}\label{2.1} 
Let $L$ be an upper semilattice and suppose $ G=\{g_{i}\mid i<\omega\}\subseteq L$. If there exist $p$, $q\in L$ such that
$$
\{g_{i}\mid i\in\omega\}\subseteq \{ x\mid x\vee{p}\geq{q}\And (\forall{y}<x)(y\vee p\not\geq q)\}
$$
then $G$ is called a Slaman-Woodin set (SW-set) for $p$, $q$ in $L$.

 If there exist $e_{0}$, $e_{1}$, $f_{0}$, $f_{1}\in L$ such that for each $ i<\omega$,
$$
g_{2i+1}=(g_{2i}\vee e_{1})\wedge f_{1}\And g_{2i+2}=(g_{2i+1}\vee e_{0})\wedge f_{0},
$$
then the function $i\mapsto g_{i}$ is called a Shore sequence for $e_{0}$, $e_{1}$, $f_{0}$, $f_{1}$ in $L$.
\end{df}

\begin{lem}\label{2.2}
Let $\mathbf a$ be a Turing degree. Let $L$ be a $\Sigma_{1}^{0}(\bo a)$-presented upper semilattice containing elements $p$, $q$, $e_{0}$, $e_{1}$, $f_{0}$, $f_{1}$, and atoms $g_{i}$ for $ i\in\omega$, such that $G=\{g_{i}\mid i<\omega\}$ is a Slaman-Woodin set for $p$, $q$ and $i\mapsto g_{i}$ is a Shore sequence for $e_{0}$, $e_{1}$, $f_{0}$, $f_{1}$. Then $\{\langle y, i\rangle\mid y=g_{i}\}\leq_T \mathbf{a}$.
\end{lem}
\begin{proof}

$$
y=g_{2i+1}\Leftrightarrow\exists x(x=g_{2i}\And  y\leq x\vee e_{1}\And  y\leq f_{1}\And  y\vee p\geq q)
$$
$$
y=g_{2i+2}\Leftrightarrow\exists x(x=g_{2i+1}\And  y\leq x\vee e_{0}\And  y\leq f_{0}\And  y\vee p\geq q)
$$
Note that the matrices of the formulas on the right hand side are positive formulas in the language with $\vee$ and $\leq$. The function $\vee$ is $\bo a$-recursive and the relation $\leq$ is $\Sigma_{1}^{0}(\mathbf{a})$. Hence the entire right hand sides are $\Sigma_{1}^{0}(\mathbf{a})$. So starting with $g_{0}$ we can find $g_n$, $\bo a$-recursively. 
\end{proof}

\begin{df}\label{2.4}
Let $a \in \mathcal{D}$. An usl $L$ is said to be of degree $\bo a$ if (1) $L$ is $\bo a$-presentable, and (2) if $\mathbf{b}\in \mathcal{D}$ and $L$ is $\mathbf{b}$-presentable then $\mathbf{a}\leq \mathbf{b}$.
\end{df}

\begin{df}\label{2.5}
Given $ U\subseteq\omega$ we define a lattice $L(U)$.

 It consists of $0,1$, atoms $\{g_{i} : i\in\omega\}$, more atoms $e_{0}$, $e_{1}$, $p$, $s$ and non-atoms $f_{0}$, $f_{1}<1$ with the properties of Lemma \ref{2.2} (taking $q=1$) and an additional element $s$ with the following property for each $n \in\omega$:
$$
n\in U\Leftrightarrow g_{n}\vee s=1.
$$
\end{df}

\begin{rem}\label{2.3} Historically, the technique of enumerating the $g_{n}$ $\bo a'$-recursively was first done in \cite{Shore82}. The idea of the improvement can be seen in \cite{Shore81} Lemma 1.11. The Slaman-Woodin conditions used to combine these ideas to get the above lemma were presented in \cite{NSS96} with a proof appearing in \cite{NSS98} Lemma 2.13(i).
 The construction of $L(U)$ was presented to the author by Slaman; see also Theorem 3.7 of \cite{NSS96}.
\end{rem}

\begin{rem} Here are some details for the proof that such a lattice $L(U)$ exists (thanks to assistance from participants in a 2008 seminar at the University of Hawai\textquoteleft i). We make $L(U)$ a height-three lattice, i.e., every element is either 0, 1, an atom or a co-atom. The atoms are $e_0$, $e_1$, $s$, and the $g_i$. The element $p$ may be either an atom or a co-atom, and is incomparable with all other elements except that $0\le p\le 1$. The co-atoms are $e_0\vee g_{2n+1}$, and $e_1\vee g_{2n}$, $f_0$, $f_1$, and $g_i\vee s$ whenever $i\not\in U$. These elements are incomparable except as forced by the above conditions. The point of including $p$ and $q$ is that $y\vee p\ge q$ is a positive statement that implies $y\not\le 0$.
\end{rem}

 The following lemma will have many applications:

\begin{lem}\label{2.6} 
Let $ U\subseteq\omega$.

\begin{enumerate}
\item[1.] $L(U)$ has degree $\mathbf{u}$.

\item[2.] If $L(U)$ is $\Sigma_{1}^{0}(\bo b)$-presentable, then $U\in\Sigma_{1}^{0}(\mathbf{b})$.

\item[3.] If $U\in\Sigma_{1}^{0}(\mathbf{b})$ then $L(U)$ is $\Sigma_{1}^{0}(\bo b)$-presentable.
\end{enumerate}

\end{lem}
\begin{proof}
1. The definition of $L(U)$ appeals to an oracle of degree $\mathbf{u}$ only and so $L(U)$ is $\bo u$-presentable. Suppose $L(U)$ is presented with degree $\mathbf{v}$. By Lemma \ref{2.2}, the relation $y=g_{i}$ is recursive in $\mathbf{v}$. Now $i\in U\leftrightarrow g_{i}\vee s\geq 1$, so since $\vee$ and $\geq$ are recursive in $\mathbf{v},\mathbf{u}\leq \mathbf{v}$.

2. We have
$$
n\in U\Leftrightarrow\exists x(x=g_{n}\And  x\vee s\geq 1)\Leftrightarrow\forall x(x=g_{n}\rightarrow x\vee s\geq 1).
$$

By Lemma \ref{2.2}, $U$ is of the form $\exists x(\triangle_{1}^{0}(\mathbf{b})\And \Sigma_{1}^{0}(\mathbf{b}))$ $U\in\Sigma_{1}^{0}(\mathbf{b})$.

3. Immediate from the fact that all clauses of the definition of $L(U)$ except ``$n\in U\leftrightarrow g_{n}\vee s\geq 1$'' are recursive.
\end{proof}

\begin{pro}\label{2.7} 
If each $\Sigma_4^0(\bo x)$-presentable bounded usl is $\Sigma_{4}^{0}(\bo y)$-presentable, then $\bo x^{(3)}\le_T \bo y^{(3)}$.
\end{pro}
\begin{proof}
Let $\mathbf{b}=\mathbf{x}^{(3)}$. Since $L(B\oplus\overline B)$ is $\Sigma^0_1(B)$-presentable, it is $\Sigma^0_1(\bo y^{(3)})$-presentable. Thus by Lemma \ref{2.6}(2), $B\oplus\overline B$ is $\Sigma^0_1(\bo y^{(3)}$ and hence $B\le_T \bo y^{(3)}$.

\end{proof}

\section*{Acknowledgments}

Thanks are due to Noam Greenberg for a correction to an early draft of this article. This material is based upon work supported by the National Science Foundation under Grants No. 0652669 and 0901020.

\end{document}